\documentclass[reqno]{amsart}
\usepackage{graphicx}
\usepackage{amssymb}
\usepackage{latexsym,amsmath}
\usepackage{hyperref}
\usepackage{amsfonts}
\usepackage{amsmath}
\usepackage{dsfont}
\begin{document}
 \title[\hfilneg ]{Katugampola Fractional Calculus With Generalized $k-$Wright
Function}

 \author[\hfil\hfilneg]{Ahmad Y. A. Salamooni, D. D. Pawar }
 \address{Ahmad Y. A. Salamooni \newline
     School of Mathematical Sciences, Swami Ramanand Teerth Marathwada University, Nanded-431606, India}
      \email{ayousss83@gmail.com}

 \address{D. D. Pawar \newline
    School of Mathematical Sciences, Swami Ramanand Teerth Marathwada University, Nanded-431606, India}
     \email{dypawar@yahoo.com}


  \begin{abstract}
 In this article, we presented some properties of the Katugampola fractional integrals
 and derivatives. Also we studied the fractional calculus properties involving Katugampola Fractional integrals
 and derivatives of generalized $k-$Wright function $_{n}\Phi_{m}^{k}(z).$\\[2mm]


  \end{abstract}

  \maketitle \numberwithin{equation}{section}
    \newtheorem{theorem}{Theorem}[section]
    \newtheorem{lemma}[theorem]{Lemma}
    \newtheorem{definition}[theorem]{Definition}
    \newtheorem{example}[theorem]{Example}

    \newtheorem{remark}[theorem]{Remark}
    \allowdisplaybreaks
 \section{Introduction and Preliminaries}
In recent years, researchers and authors have introduced a new fractional  integrators  operators
and fractional derivatives operators which are generalizations of the famous
 Riemann-Leuville and  the Hadamard-type, for more details see [1-5] and references therein.
\\\textbf{\ Definition 1.[5]} Let $\Omega=[a,b],$ the Katugampola fractional integrals $_{\rho}I_{0+}^{\gamma}\varphi~and~_{\rho}I_{-}^{\gamma}\varphi$ of order $\gamma\in\mathbb{C}(\mathfrak{R}(\gamma)>0)$
are defined for $\rho>0,~a=0~and~b=\infty$ as
\begin{align}
(_{\rho}I_{0+}^{\gamma}\varphi)(s)=\frac{\rho^{1-\gamma}}{\Gamma(\gamma)}
\int_{0}^{s}\frac{\tau^{\rho-1}\varphi(\tau)}{(s^{\rho}-\tau^{\rho})^{1-\gamma}}d\tau,\quad(s>0),\label{e1.1}
\end{align}
and
\begin{align}
(_{\rho}I_{-}^{\gamma}\varphi)(s)=\frac{\rho^{1-\gamma}}{\Gamma(\gamma)}
\int_{s}^{\infty}\frac{\tau^{\rho-1}\varphi(\tau)}{(\tau^{\rho}-s^{\rho})^{1-\gamma}}d\tau,\quad(s>0),\label{e1.2}
\end{align}
and the corresponding Katugampola fractional derivatives $_{\rho}D_{0+}^{\gamma}\varphi~and~_{\rho}D_{-}^{\gamma}\varphi$
are defined with $\big(n=1+[\mathfrak{R}(\gamma)]\big)$ as
\begin{align}
(_{\rho}D_{0+}^{\gamma}\varphi)(s)&:=\big(s^{1-\rho}\frac{d}{ds}\big)^{1+[\mathfrak{R}(\gamma)]}
\big(_{\rho}I_{0+}^{1-\gamma+[\mathfrak{R}(\gamma)]}\varphi\big)(s)\nonumber\\&
=\frac{\rho^{\gamma-[\mathfrak{R}(\gamma)]}}{\Gamma(1-\gamma+[\mathfrak{R}(\gamma)])}\big(s^{1-\rho}\frac{d}{ds}\big)^{1+[\mathfrak{R}(\gamma)]}
\int_{0}^{s}\frac{\tau^{\rho-1}\varphi(\tau)}{(s^{\rho}-\tau^{\rho})^{\gamma-[\mathfrak{R}(\gamma)]}}d\tau,\quad(s>0),\label{e1.3}
\end{align}
and
\begin{align}
(_{\rho}D_{-}^{\gamma}\varphi)(s)&:=\big(-s^{1-\rho}\frac{d}{ds}\big)^{1+[\mathfrak{R}(\gamma)]}
\big(_{\rho}I_{-}^{1-\gamma+[\mathfrak{R}(\gamma)]}\varphi\big)(s)\nonumber\\&
=\frac{\rho^{\gamma-[\mathfrak{R}(\gamma)]}}{\Gamma(1-\gamma+[\mathfrak{R}(\gamma)])}
\big(-s^{1-\rho}\frac{d}{ds}\big)^{1+[\mathfrak{R}(\gamma)]}
\int_{s}^{\infty}\frac{\tau^{\rho-1}\varphi(\tau)}{(\tau^{\rho}-s^{\rho})^{\gamma-[\mathfrak{R}(\gamma)]}}d\tau,\quad(s>0).\label{e1.4}
\end{align}
\\
\textbf{\ Definition 2.[6]}  The generalized $K-$Gamma function $\Gamma_{k}(y)$ defined by
\begin{align}
\Gamma_{k}(y)=\lim_{n\rightarrow\infty}\frac{n!k^{n}(nk)^{\frac{y}{k}-1}}{(y)_{n,k}},
\quad(k>0;~y\in\mathbb{C}\setminus k\mathbb{Z^{-}}),\label{e1.5}
\end{align}
where $(y)_{n,k}$ is the $k-$Pochhammer symbol given as
\begin{align}
(y)_{n,k}:=\left\{\begin{matrix}\frac{\Gamma_{k}(y+nk)}{\Gamma_{k}(y)}
\quad\quad\quad\quad\quad\quad\quad\quad\quad(k\in\mathbb{R};~y\in\mathbb{C}\setminus \{0\})\\\\
y(y+k)(y+2k)...(y+(n-1)k)\quad\quad(n\in\mathbb{N^{+}};~y\in\mathbb{C})\end{matrix}\right.\label{e1.6}
\end{align}
and for $\mathfrak{R}(y)>0,$ the $K-$Gamma function $\Gamma_{k}(y)$ defined by the integral
\begin{align}
\Gamma_{k}(y)=\int_{0}^{\infty}x^{y-1}e^{-\frac{x^{k}}{k}}dx\quad\quad\quad\quad\quad\quad\quad\quad\label{e1.7}
\end{align}
this given relation with Euler's Gamma function as
\begin{align}
\Gamma_{k}(y)=k^{\frac{y}{k}-1}\Gamma(\frac{y}{k}).\quad\quad\quad\quad\quad\quad\quad\quad\label{e1.8}
\end{align}
Also [7],
\begin{align}
\Gamma(1-y)\Gamma(y)=\frac{\pi}{\sin(y\pi)}.\quad\quad\quad\quad\quad\quad\quad\quad\label{e1.9}
\end{align}
\\
\textbf{\ Definition 3. [8]} The Beta function $B(\upsilon,\omega)$ is defined as
\begin{align}\nonumber
B(\upsilon,\omega)&= \int_{0}^{1}z^{\upsilon-1}(1-z)^{\omega-1}dz,\quad\mathfrak{R}(\upsilon)>0,\quad\mathfrak{R}(\omega)>0,\\&
=\frac{\Gamma(\upsilon)\Gamma(\omega)}{\Gamma(\upsilon+\omega)}\label{e1.10}
\end{align}
Furthermore,
\begin{align}\nonumber
&\int_{\hat{x}}^{\infty}(z-\hat{x})^{\upsilon-1}(z-\hat{y})^{\omega-1}dz=
(\hat{x}-\hat{y})^{\upsilon+\omega-1} B(\upsilon,1-\upsilon-\omega), \\&\quad\quad\quad\quad\quad\quad \hat{x}>\hat{y},\quad0<\mathfrak{R}(\upsilon)<1-\mathfrak{R}(\omega).\label{e1.11}
\end{align}
\\
\par Recently the Generalized $K-$Wright function introduced by (Gehlot and Prajapati [9]) which is defined as following:
\\\textbf{\ Definition 4.}  For $k\in\mathbb{R^{+}};~z\in\mathbb{C};~p_{i},q_{j}\in\mathbb{C},
~\alpha_{i},\beta_{j}\in\mathbb{R}~(\alpha_{i},\beta_{j}\neq0;~i=1,2,...,n;~j=1,2,...,m)$ and
 $(p_{i}+\alpha_{i}r),~(q_{j}+\beta_{j}r)\in\mathbb{C}\setminus k\mathbb{Z^{-}},$ the
 generalized $k-$Wright function $_{n}\Phi_{m}^{k}$ is defined by
\begin{align}
_{n}\Phi_{m}^{k}(z)=~_{n}\Phi_{m}^{k}\Bigg[\begin{matrix}(p_{i},\alpha_{i})_{1,n}\\(q_{j},\beta_{j})_{1,m}
\end{matrix}\Big|z\Bigg]=\sum_{r=0}^{\infty}\frac{\prod_{i=1}^{n}~\Gamma_{k}(p_{i}+\alpha_{i}r)}{\prod_{j=1}^{m}\Gamma_{k}(q_{j}+
\beta_{j}r)}\frac{z^{r}}{r!},\label{e1.12}
\end{align}
with the convergence conditions describing as
\begin{align}\nonumber
\Delta=\sum_{j=1}^{m}\big(\frac{\beta_{j}}{k}\big)-\sum_{i=1}^{n}
\big(\frac{\alpha_{i}}{k}\big);\mu=\prod_{i=1}^{n}\big|\frac{\alpha_{i}}{k}\big|^{-\frac{\alpha_{i}}{k}}
\prod_{j=1}^{m}\big|\frac{\beta_{j}}{k}\big|^{\frac{\beta_{j}}{k}};\nu=
\sum_{j=1}^{m}\big(\frac{q_{j}}{k}\big)-\sum_{i=1}^{n}\big(\frac{p_{i}}{k}\big)+\frac{n-m}{2}
\end{align}
\\\textbf{\ Lemma 1. [9]}
For $k\in\mathbb{R^{+}};~z\in\mathbb{C};~p_{i},q_{j}\in\mathbb{C},~\alpha_{i},
\beta_{j}\in\mathbb{R}~(\alpha_{i},\beta_{j}\neq0;~i=1,2,...,n;~j=1,2,...,m)$ and
$(p_{i}+\alpha_{i}r),~(q_{j}+\beta_{j}r)\in\mathbb{C}\setminus k\mathbb{Z^{-}}$
\par $(1)$ If $\Delta>-1,$ then series $(1.12)$ is absolutely convergent for all $z \in\mathbb{C}$ and
\par$~~~\quad$ generalized $k-$Wright function $_{n}\Phi_{m}^{k}(z)$ is an entire function of $z.$
\par $(2)$ If $\Delta=-1,$ then series $(1.12)$ is absolutely convergent for all $|z| < \mu$ and of $$|z| = \mu,\mathfrak{R}(\mu)>\frac{1}{2}.$$

  \section{Properties of Katugampola Fractional integral and derivative}
\par In this section,  we investigated some properties of the Katugampola fractional integrals and derivatives $(1.1),(1.2)~and~(1.3),(1.4)$ for the power function $\varphi(s)=s^{\alpha-1}$ and the exponential function $e^{-\lambda~s^{\rho}}.$
\\\textbf{\ Lemma 2.} Let $\rho>0,\mathfrak{R}(\gamma)\geqq0~~ and~~ n=1+[\mathfrak{R}(\gamma)]$
\par $(1)~~If~\mathfrak{R}(\alpha)>0,$ then
\begin{align}
(_{\rho}I_{0+}^{\gamma}\tau^{\alpha-1})(s)=
\frac{\rho^{-\gamma}\Gamma(1+\frac{\alpha-1}{\rho})}{\Gamma(1+\frac{\alpha-1}{\rho}+\gamma)}s^{\rho\gamma+(\alpha-1)}
\quad(\mathfrak{R}(\gamma)\geq0;~\mathfrak{R}(\alpha)>0)\label{e2.1}
\end{align}
\begin{align}
\quad\quad(_{\rho}D_{0+}^{\gamma}\tau^{\alpha-1})(s)=
\frac{\rho^{\gamma-n}\Gamma(1+\frac{\alpha-1}{\rho})}{\Gamma(1+\frac{\alpha-1}{\rho}-\gamma)}s^{(\alpha-1)-\rho\gamma}
\quad(\mathfrak{R}(\gamma)\geqq0;~\mathfrak{R}(\alpha)>0).\label{e2.2}
\end{align}
\par $(2)~~If~\alpha\in\mathbb{C},$ then
\begin{align}
\quad\quad(_{\rho}I_{-}^{\gamma}\tau^{\alpha-1})(s)=
\frac{\rho^{-\gamma}\Gamma(\frac{1-\alpha}{\rho}-\gamma)}{\Gamma(\frac{1-\alpha}{\rho})}s^{\rho\gamma+(\alpha-1)}
\quad(\mathfrak{R}(\gamma)\geq0;~\mathfrak{R}(\gamma+\alpha)<1)\label{e2.3}
\end{align}
\begin{align}
\quad\quad\quad(_{\rho}D_{-}^{\gamma}\tau^{\alpha-1})(s)=
\frac{\rho^{\gamma-n}\Gamma(\frac{1-\alpha}{\rho}+\gamma)}{\Gamma(\frac{1-\alpha}{\rho})}s^{(\alpha-1)-\rho\gamma}
\quad(\mathfrak{R}(\gamma)\geqq0;~\mathfrak{R}(\gamma+\alpha-[\mathfrak{R}(\gamma)])<1).\label{e2.4}
\end{align}
\par $(3)~~If~\mathfrak{R}(\lambda)>0,$ then
\begin{align}
(_{\rho}I_{-}^{\gamma}e^{-\lambda\tau^{\rho}})(s)=(\lambda\rho)^{-\gamma}e^{-\lambda~s^{\rho}}
\quad\quad(\mathfrak{R}(\gamma)\geq0)\label{e2.5}
\end{align}
\begin{align}
(_{\rho}D_{-}^{\gamma}e^{-\lambda\tau^{\rho}})(s)=(\lambda\rho)^{\gamma}e^{-\lambda~s^{\rho}}
\quad\quad(\mathfrak{R}(\gamma)\geqq0).\label{e2.6}
\end{align}
\\\textbf{\ Proof.} To prove this Lemma, let the substitution $x=\frac{\tau^{\rho}}{s^{\rho}},$ in parts (1) and (2).
\par $(1)$ Firstly, by the equation (1.1) and the given substitution we have
\begin{align}\nonumber
(_{\rho}I_{0+}^{\gamma}\tau^{\alpha-1})(s)&=\frac{\rho^{-\gamma}s^{\rho\gamma+\alpha-1}}{\Gamma(\gamma)}
\int_{0}^{1}\frac{x^{\frac{\alpha-1}{\rho}}}{(1-x)^{1-\gamma}}dx\\&\nonumber
=\frac{\rho^{-\gamma}s^{\rho\gamma+\alpha-1}}{\Gamma(\gamma)}~B\big(\gamma,1+\frac{\alpha-1}{\rho}\big),
\end{align}
now, using equation (1.10), we obtain the result (2.1).
$\quad\quad\quad\quad\quad\quad\quad\quad\quad\quad\quad\quad\quad\quad\Box$
\\
Secondly, by the equation (1.3), the given substitution and by  using the result (2.1), we have
\begin{align}\nonumber
\quad\quad(_{\rho}D_{0+}^{\gamma}\tau^{\alpha-1})(s)&=
\big(s^{1-\rho}\frac{d}{ds}\big)^{n}\big(_{\rho}I_{0+}^{n-\gamma}\tau^{\alpha-1}\big)(s)
\\&\nonumber
=\frac{\rho^{\gamma-n}\Gamma(1+\frac{\alpha-1}{\rho})}{\Gamma(1+\frac{\alpha-1}{\rho}+n-\gamma)}
\big(s^{1-\rho}\frac{d}{ds}\big)^{n}s^{\rho(n-\gamma)+\alpha-1}
\\&\nonumber
=\frac{\rho^{\gamma-n}\Gamma(1+\frac{\alpha-1}{\rho})}{\Gamma(1+\frac{\alpha-1}{\rho}-\gamma)}s^{(\alpha-1)-\rho\gamma}\quad\quad\quad
\quad\quad\quad\quad\quad\quad\quad\quad\quad\quad\quad\Box
\end{align}
 \par $(2)$ Firstly, by the equation (1.2) and the given substitution we have
\begin{align}\nonumber
\quad\quad\quad(_{\rho}I_{-}^{\gamma}\tau^{\alpha-1})(s)=\frac{\rho^{-\gamma}s^{\rho\gamma+\alpha-1}}{\Gamma(\gamma)}
\int_{1}^{\infty}x^{\frac{\alpha-1}{\rho}}(x-1)^{\gamma-1}dx,
\end{align}
now, using the equation (1.11), with $\hat{x}=1~~and~~\hat{y}=0,$ we obtain
\begin{align}\nonumber
(_{\rho}I_{-}^{\gamma}\tau^{\alpha-1})(s)=\frac{\rho^{-\gamma}s^{\rho\gamma+\alpha-1}}{\Gamma(\gamma)}
B\big(\gamma,1-\gamma-(1+\frac{\alpha-1}{\rho})\big),
\end{align}
by using equation (1.10), we obtain the result (2.3).
$\quad\quad\quad\quad\quad\quad\quad\quad\quad\quad\quad\quad\quad\quad\Box$
\\\\
Secondly, by the equation (1.4), the given substitution and by using the result (2.3), we have
\begin{align}\nonumber
(_{\rho}D_{-}^{\gamma}\tau^{\alpha-1})(s)&=
\big(-s^{1-\rho}\frac{d}{ds}\big)^{n}\big(_{\rho}I_{-}^{n-\gamma}\tau^{\alpha-1}\big)(s)
\\&\nonumber
=\frac{(-1)^{n}\rho^{\gamma-n}\Gamma(\frac{1-\alpha}{\rho}+\gamma-n)}
{\Gamma(\frac{1-\alpha}{\rho})}\big(s^{1-\rho}\frac{d}{ds}\big)^{n}
s^{\rho(n-\gamma)+\alpha-1}
\\&
=\frac{(-1)^{n}\rho^{\gamma-n}}{\Gamma(\frac{1-\alpha}{\rho})}\frac{\Gamma(\frac{1-\alpha}{\rho}+\gamma-n)
\Gamma(1-[\frac{1-\alpha}{\rho}+\gamma-n])}{\Gamma(1-[\gamma-\frac{\alpha-1}{\rho}])}\label{e2.7}
\end{align}
Also, by using (1.9), we have
\begin{align}
\quad&\Gamma(\frac{1-\alpha}{\rho}+\gamma-n)\Gamma(1-[\frac{1-\alpha}{\rho}+\gamma-n])=
\frac{\pi}{\sin([\frac{1-\alpha}{\rho}+\gamma-n]\pi)}=
\frac{(-1)^{n}\pi}{\sin([\gamma-\frac{\alpha-1}{\rho}]\pi)}\label{e2.8}
\\&\nonumber and
\\&\frac{1}{\Gamma(1-[\gamma-\frac{\alpha-1}{\rho}])}=
\frac{\Gamma(\gamma-\frac{\alpha-1}{\rho})}{\Gamma(\gamma-\frac{\alpha-1}{\rho})
\Gamma(1-[\gamma-\frac{\alpha-1}{\rho}])}=
\frac{\Gamma(\gamma-\frac{\alpha-1}{\rho})}{\pi}\sin([\gamma-\frac{\alpha-1}{\rho}]\pi)\label{e2.9}
\end{align}
Substituting relations (2.8) and (2.9) in (2.7), we obtain (2.4). $\quad\quad\quad\quad\quad\quad\quad\quad\quad\quad\quad\Box$
\par $(3)$ For this part let the substitution $x=\tau^{\rho}-s^{\rho}.$
\\\\Firstly, by the equation (1.2) and the given substitution in this part we have
\begin{align}\nonumber
(_{\rho}I_{-}^{\gamma}e^{-\lambda\tau^{\rho}})(s)=
\frac{\rho^{-\gamma}}{\Gamma(\gamma)}e^{-\lambda~s^{\rho}}\int_{0}^{\infty}e^{-\lambda~x}x^{\gamma-1}dx,
\end{align}
then use the substitution $\vartheta=\lambda~x,$ we obtain
\begin{align}\nonumber
(_{\rho}I_{-}^{\gamma}e^{-\lambda\tau^{\rho}})(s)=\frac{\rho^{-\gamma}}{\Gamma(\gamma)}e^{-\lambda~s^{\rho}}\lambda^{-\gamma}
\int_{0}^{\infty}e^{-\vartheta}\vartheta^{\gamma-1}d\vartheta,
\end{align}
since $~~\int_{0}^{\infty}e^{-\vartheta}\vartheta^{\gamma-1}d\vartheta=\Gamma(\gamma)~~~[7],$ then the result is satisfy.$\quad\quad\quad\quad\quad\quad\quad\quad\quad\quad\quad\quad\quad\Box$\\\\
Secondly, by the equation (1.4) and by using the result (2.5), we have
\begin{align}\nonumber
(_{\rho}D_{-}^{\gamma}e^{-\lambda\tau^{\rho}})(s)&=
\big(-s^{1-\rho}\frac{d}{ds}\big)^{n}\big(_{\rho}I_{-}^{n-\gamma}e^{-\lambda\tau^{\rho}}\big)(s)
\\&\nonumber
=(-1)^{n}\big(s^{1-\rho}\frac{d}{ds}\big)^{n}\big((\lambda\rho)^{\gamma-n}e^{-\lambda~s^{\rho}}\big)
\\&\nonumber
=(-1)^{n}~s^{(1-\rho)n}~(\lambda\rho)^{\gamma-n}\big(\frac{d^{n}}{ds^{n}}e^{-\lambda~s^{\rho}}\big)
\\&\nonumber
=(\lambda\rho)^{\gamma}e^{-\lambda~s^{\rho}}\quad\quad\quad\quad\quad\quad
\quad\quad\quad\quad\quad\quad\quad\quad\quad\quad\quad\quad\quad\quad\quad\Box
\end{align}
\\\textbf{\ Remark 1.} $(a)$ In Lemma 2, if the power function is $\varphi(s)=\big(\frac{s^{\rho}}{\rho}\big)^{\alpha-1},$ then
\par $(1)~~If~\mathfrak{R}(\alpha)>0,$ then
\begin{align}\nonumber
\bigg(~_{\rho}I_{0+}^{\gamma}\big(\frac{\tau^{\rho}}{\rho}\big)^{\alpha-1}\bigg)(s)=
\frac{\Gamma(\alpha)}{\Gamma(\alpha+\gamma)}\big(\frac{s^{\rho}}{\rho}\big)^{\alpha+\gamma-1}
\quad(\mathfrak{R}(\gamma)\geq0;~\mathfrak{R}(\alpha)>0)
\end{align}
\begin{align}\nonumber
\quad\quad\bigg(~_{\rho}D_{0+}^{\gamma}\big(\frac{\tau^{\rho}}{\rho}\big)^{\alpha-1}\bigg)(s)=
\frac{\Gamma(\alpha)}{\Gamma(\alpha-\gamma)}\big(\frac{s^{\rho}}{\rho}\big)^{\alpha-\gamma-1}
\quad(\mathfrak{R}(\gamma)\geqq0;~\mathfrak{R}(\alpha)>0).
\end{align}
\par $(2)~~If~\alpha\in\mathbb{C},$ then
\begin{align}\nonumber
\quad\quad\bigg(~_{\rho}I_{-}^{\gamma}\big(\frac{\tau^{\rho}}{\rho}\big)^{\alpha-1}\bigg)(s)=
\frac{\Gamma(1-\gamma-\alpha)}{\Gamma(1-\alpha)}\big(\frac{s^{\rho}}{\rho}\big)^{\alpha+\gamma-1}
\quad(\mathfrak{R}(\gamma)\geq0;~\mathfrak{R}(\gamma+\alpha)<1)
\end{align}
\begin{align}\nonumber
\quad\quad\quad\bigg(~_{\rho}D_{-}^{\gamma}\big(\frac{\tau^{\rho}}{\rho}\big)^{\alpha-1}\bigg)(s)=
\frac{\Gamma(1+\gamma-\alpha)}{\Gamma(1-\alpha)}\big(\frac{s^{\rho}}{\rho}\big)^{\alpha-\gamma-1}
\quad(\mathfrak{R}(\gamma)\geqq0;~\mathfrak{R}(\gamma+\alpha-[\mathfrak{R}(\gamma)])<1).
\end{align}
\\
$(b)$ If $\mathfrak{R}(\alpha)>\mathfrak{R}(\gamma)>0,$ then
\begin{align}
(_{\rho}I_{-}^{\gamma}\tau^{-\alpha})(s)=\frac{\rho^{-\gamma}\Gamma(\frac{\alpha}{\rho}-\gamma)}{\Gamma(\frac{\alpha}{\rho})}s^{\rho\gamma-\alpha}
\end{align}
\section{Katugampola Fractional integration  for Generalized $k-$Wright
Function}
\par In this section, we established the Katugampola fractional integration for generalized $k-$Wright
function (1.12).
\\\textbf{\ Theorem 1.}
Let $\gamma,~ \alpha\in\mathbb{C}$ such that $\mathfrak{R}(\gamma)>0,~ \mathfrak{R}(\alpha)>0;~ \lambda\in\mathbb{C},~ \rho>0,~\nu>0,$
then for $\Delta>-1,$ the Katugampola fractional integration $_{\rho}I_{0+}^{\gamma}$ for generalized $k-$Wright
function $_{n}\Phi_{m}^{k}(z)$ is given as
\begin{align}\nonumber
\Bigg(~_{\rho}I_{0+}^{\gamma}&\Bigg(\tau^{\frac{\alpha}{k}-1}~_{n}\Phi_{m}^{k}
\Bigg[\begin{matrix}(p_{i},\alpha_{i})_{1,n}\\(q_{j},\beta_{j})_{1,m}
\end{matrix}\Big|~\lambda~\tau^{\frac{\nu}{k}}\Bigg]\Bigg) \Bigg)(s)\\&
=(\frac{k}{\rho})^{\gamma}~s^{\frac{\alpha}{k}+\rho\gamma-1}
~_{n+1}\Phi_{m+1}^{k}\Bigg[\begin{matrix}\big(p_{i},\alpha_{i}\big)_{1,n},~\big(\frac{1}{\rho}(\alpha+(\rho-1)k),\frac{\nu}{\rho}\big)\quad \\
\big(q_{j},\beta_{j}\big)_{1,m},~\big(\frac{1}{\rho}(\alpha+(\rho(\gamma+1)-1)k),\frac{\nu}{\rho}\big)
\end{matrix}\Bigg|~\lambda~s^{\frac{\nu}{k}}\Bigg]\label{e3.1}
\end{align}
\\\textbf{\ Proof.} According to Lemma 1, a generalized $k-$Wright
functions in both side of the equation (3.1), exist for $s>0.$
We consider that
\begin{align}\nonumber
M\equiv\Bigg(~_{\rho}I_{0+}^{\gamma}\Bigg(\tau^{\frac{\alpha}{k}-1}~_{n}\Phi_{m}^{k}
\Bigg[\begin{matrix}(p_{i},\alpha_{i})_{1,n}\\(q_{j},\beta_{j})_{1,m}
\end{matrix}\Big|~\lambda~\tau^{\frac{\nu}{k}}\Bigg]\Bigg) \Bigg)(s)
\end{align}
using (1.12), we can write the above equation as
\begin{align}\nonumber
M\equiv\Bigg(~_{\rho}I_{0+}^{\gamma}\Bigg(\tau^{\frac{\alpha}{k}-1}~\sum_{r=0}^{\infty}
\frac{\prod_{i=1}^{n}~\Gamma_{k}(p_{i}+\alpha_{i}r)}{\prod_{j=1}^{m}\Gamma_{k}(q_{j}+
\beta_{j}r)}\frac{(\lambda~\tau^{\frac{\nu}{k}})^{r}}{r!}\Bigg) \Bigg)(s)
\end{align}
now, using the integration of the series term- by term we obtain
\begin{align}\nonumber
M\equiv~\sum_{r=0}^{\infty}
\frac{\prod_{i=1}^{n}~\Gamma_{k}(p_{i}+\alpha_{i}r)}{\prod_{j=1}^{m}\Gamma_{k}(q_{j}+
\beta_{j}r)}\frac{(\lambda)^{r}}{r!}\Big(~_{\rho}I_{0+}^{\gamma}\Big(\tau^{\frac{\alpha}{k}+\frac{\nu r}{k}-1}\Big) \Big)(s)
\end{align}
applying (2.1), the above equation reduces to,
\begin{align}\nonumber
M\equiv~\sum_{r=0}^{\infty}
\frac{\prod_{i=1}^{n}~\Gamma_{k}(p_{i}+\alpha_{i}r)}{\prod_{j=1}^{m}\Gamma_{k}(q_{j}+
\beta_{j}r)}\frac{(\lambda)^{r}}{r!}\frac{\rho^{-\gamma}\Gamma(1+\frac{\frac{\alpha}{k}+
\frac{\nu r}{k}-1}{\rho})}{\Gamma(1+\frac{\frac{\alpha}{k}+\frac{\nu r}{k}-1}{\rho}+\gamma)}
s^{\frac{\alpha+\nu r}{k}+\rho\gamma-1}
\end{align}
using (1.8), we obtain
\begin{align}\nonumber
M\equiv~(\frac{k}{\rho})^{\gamma}~s^{\frac{\alpha}{k}+\rho\gamma-1}
~_{n+1}\Phi_{m+1}^{k}\Bigg[\begin{matrix}\big(p_{i},\alpha_{i}\big)_{1,n},~\big(\frac{1}{\rho}(\alpha+(\rho-1)k),\frac{\nu}{\rho}\big)\quad \\
\big(q_{j},\beta_{j}\big)_{1,m},~\big(\frac{1}{\rho}(\alpha+(\rho(\gamma+1)-1)k),\frac{\nu}{\rho}\big)
\end{matrix}\Bigg|~\lambda~s^{\frac{\nu}{k}}\Bigg]\quad\quad\quad\Box
\end{align}
\\\textbf{\ Theorem 2.}
Let $\gamma,~ \alpha\in\mathbb{C}$ such that $\mathfrak{R}(\gamma)>0,~ \mathfrak{R}(\alpha)>0;~ \lambda\in\mathbb{C},~ \rho>0,~\nu>0,$
then for $\Delta>-1,$ the Katugampola fractional integration $_{\rho}I_{-}^{\gamma}$ for generalized $k-$Wright
function $_{n}\Phi_{m}^{k}(z)$ is given as
\begin{align}\nonumber
\Bigg(~_{\rho}I_{-}^{\gamma}&\Bigg(\tau^{-\frac{\alpha}{k}}~_{n}\Phi_{m}^{k}
\Bigg[\begin{matrix}(p_{i},\alpha_{i})_{1,n}\\(q_{j},\beta_{j})_{1,m}
\end{matrix}\Big|~\lambda~\tau^{-\frac{\nu}{k}}\Bigg]\Bigg) \Bigg)(s)\\&
=(\frac{k}{\rho})^{\gamma}~s^{\rho\gamma-\frac{\alpha}{k}}
~_{n+1}\Phi_{m+1}^{k}\Bigg[\begin{matrix}\big(p_{i},\alpha_{i}\big)_{1,n},~\big(\frac{\alpha}{\rho}-k\gamma,\frac{\nu}{\rho}\big)\quad \\
\big(q_{j},\beta_{j}\big)_{1,m},~\big(\frac{\alpha}{\rho},\frac{\nu}{\rho}\big)
\end{matrix}\Bigg|~\lambda~s^{-\frac{\nu}{k}}\Bigg]\label{e3.2}
\end{align}
\\\textbf{\ Proof.} According to Lemma 1, a generalized $k-$Wright
functions in both side of the equation (3.2), exist for $s>0.$
We consider that
\begin{align}\nonumber
N\equiv\Bigg(~_{\rho}I_{-}^{\gamma}\Bigg(\tau^{-\frac{\alpha}{k}}~_{n}\Phi_{m}^{k}
\Bigg[\begin{matrix}(p_{i},\alpha_{i})_{1,n}\\(q_{j},\beta_{j})_{1,m}
\end{matrix}\Big|~\lambda~\tau^{-\frac{\nu}{k}}\Bigg]\Bigg) \Bigg)(s)
\end{align}
using (1.12), we can write the above equation as
\begin{align}\nonumber
N\equiv~\sum_{r=0}^{\infty}
\frac{\prod_{i=1}^{n}~\Gamma_{k}(p_{i}+\alpha_{i}r)}{\prod_{j=1}^{m}\Gamma_{k}(q_{j}+
\beta_{j}r)}\frac{(\lambda)^{r}}{r!}\Big(~_{\rho}I_{-}^{\gamma}\Big(\tau^{-\frac{\alpha+\nu r}{k}}\Big) \Big)(s)
\end{align}
applying (2.10), the above equation reduces to,
\begin{align}\nonumber
N\equiv~\sum_{r=0}^{\infty}
\frac{\prod_{i=1}^{n}~\Gamma_{k}(p_{i}+\alpha_{i}r)}{\prod_{j=1}^{m}\Gamma_{k}(q_{j}+
\beta_{j}r)}\frac{(\lambda)^{r}}{r!}\frac{\rho^{-\gamma}\Gamma(\frac{\frac{\alpha+\nu r}{k}}
{\rho}-\gamma)}{\Gamma(\frac{\frac{\alpha+\nu r}{k}}
{\rho})}
s^{\rho\gamma-\frac{\alpha+\nu r}{k}}
\end{align}
using (1.8), we obtain
\begin{align}\nonumber
\quad\quad~N\equiv~(\frac{k}{\rho})^{\gamma}~s^{\rho\gamma-\frac{\alpha}{k}}
~_{n+1}\Phi_{m+1}^{k}\Bigg[\begin{matrix}\big(p_{i},\alpha_{i}\big)_{1,n},~\big(\frac{\alpha}{\rho}-k\gamma,\frac{\nu}{\rho}\big)\quad \\
\big(q_{j},\beta_{j}\big)_{1,m},~\big(\frac{\alpha}{\rho},\frac{\nu}{\rho}\big)
\end{matrix}\Bigg|~\lambda~s^{-\frac{\nu}{k}}\Bigg]\quad\quad\quad\quad\quad\quad\quad\quad\Box
\end{align}
\section{Katugampola Fractional differentiation  for Generalized $k-$Wright
Function}
\par This section deals with the Katugampola fractional differentiation for generalized $k-$Wright
function (1.12).
\\\textbf{\ Theorem 3.}
Let $\gamma,~ \alpha\in\mathbb{C}$ such that $\mathfrak{R}(\gamma)>0,~ \mathfrak{R}(\alpha)>0;~ \lambda\in\mathbb{C},~ \rho>0,~\nu>0,$
then for $\Delta>-1,$ the Katugampola fractional integration $_{\rho}D_{0+}^{\gamma}$ for generalized $k-$Wright
function $_{n}\Phi_{m}^{k}(z)$ is given as
\begin{align}\nonumber
\Bigg(~_{\rho}D_{0+}^{\gamma}&\Bigg(\tau^{\frac{\alpha}{k}-1}~_{n}\Phi_{m}^{k}
\Bigg[\begin{matrix}(p_{i},\alpha_{i})_{1,n}\\(q_{j},\beta_{j})_{1,m}
\end{matrix}\Big|~\lambda~\tau^{\frac{\nu}{k}}\Bigg]\Bigg) \Bigg)(s)\\&
=(\frac{k}{\rho})^{-\gamma}~s^{\frac{\alpha}{k}-\rho\gamma-1}
~_{n+1}\Phi_{m+1}^{k}\Bigg[\begin{matrix}\big(p_{i},\alpha_{i}\big)_{1,n},~\big(\frac{1}{\rho}(\alpha+(\rho-1)k),\frac{\nu}{\rho}\big)\quad \\
\big(q_{j},\beta_{j}\big)_{1,m},~\big(\frac{1}{\rho}(\alpha+(\rho(1-\gamma)-1)k),\frac{\nu}{\rho}\big)
\end{matrix}\Bigg|~\lambda~s^{\frac{\nu}{k}}\Bigg]\label{e4.1}
\end{align}
\\\textbf{\ Proof.} According to Lemma 1, a generalized $k-$Wright
functions in both side of the equation (4.1), exist for $s>0.$
Let $n=1+[\mathfrak{R}(\gamma)],$ then we consider that
\begin{align}\nonumber
P\equiv\Bigg(~_{\rho}D_{0+}^{\gamma}\Bigg(\tau^{\frac{\alpha}{k}-1}~_{n}\Phi_{m}^{k}
\Bigg[\begin{matrix}(p_{i},\alpha_{i})_{1,n}\\(q_{j},\beta_{j})_{1,m}
\end{matrix}\Big|~\lambda~\tau^{\frac{\nu}{k}}\Bigg]\Bigg) \Bigg)(s)
\end{align}
using (1.3), we have
\begin{align}\nonumber
P\equiv\big(s^{1-\rho}\frac{d}{ds}\big)^{n}\Bigg(~_{\rho}I_{0+}^{n-\gamma}\Bigg(\tau^{\frac{\alpha}{k}-1}~_{n}\Phi_{m}^{k}
\Bigg[\begin{matrix}(p_{i},\alpha_{i})_{1,n}\\(q_{j},\beta_{j})_{1,m}
\end{matrix}\Big|~\lambda~\tau^{\frac{\nu}{k}}\Bigg]\Bigg) \Bigg)(s)
\end{align}
using Theorem 1, we obtain
\begin{align}\nonumber
P\equiv\big(s^{1-\rho}\frac{d}{ds}\big)^{n}\Bigg((\frac{k}{\rho})^{n-\gamma}~s^{\frac{\alpha}{k}+\rho(n-\gamma)-1}
~_{n+1}\Phi_{m+1}^{k}\Bigg[\begin{matrix}\big(p_{i},\alpha_{i}\big)_{1,n},~\big(\frac{1}{\rho}(\alpha+(\rho-1)k),\frac{\nu}{\rho}\big)\quad \\
\big(q_{j},\beta_{j}\big)_{1,m},~\big(\frac{1}{\rho}(\alpha+(\rho(n-\gamma+1)-1)k),\frac{\nu}{\rho}\big)
\end{matrix}\Bigg|~\lambda~s^{\frac{\nu}{k}}\Bigg]\Bigg)
\end{align}
using (1.12), we can write the above equation as
\begin{align}\nonumber
P\equiv~(\frac{k}{\rho})^{n-\gamma}~\sum_{r=0}^{\infty}
\frac{\prod_{i=1}^{n}~\Gamma_{k}(p_{i}+\alpha_{i}r)\Gamma_{k}(\frac{1}{\rho}(\alpha+(\rho-1)k)+\frac{\nu}{\rho}r)}
{\prod_{j=1}^{m}~\Gamma_{k}(q_{j}+\beta_{j}r)\Gamma_{k}(\frac{1}{\rho}(\alpha+(\rho(n-\gamma+1)-1)k)+\frac{\nu}{\rho}r)}
\frac{(\lambda)^{r}}{r!}\big(s^{1-\rho}\frac{d}{ds}\big)^{n}\big(~s^{\frac{\alpha}{k}+\frac{\nu}{k}+\rho(n-\gamma)-1}\big)
\end{align}
we can write the above equation as
\begin{align}\nonumber
P\equiv~k^{n-\gamma}~\rho^{\gamma}~\sum_{r=0}^{\infty}&
\frac{\prod_{i=1}^{n}~\Gamma_{k}(p_{i}+\alpha_{i}r)\Gamma_{k}(\frac{1}{\rho}(\alpha+(\rho-1)k)+\frac{\nu}{\rho}r)}
{\prod_{j=1}^{m}~\Gamma_{k}(q_{j}+\beta_{j}r)\Gamma_{k}(\frac{1}{\rho}(\alpha+(\rho(n-\gamma+1)-1)k)+\frac{\nu}{\rho}r)}
\frac{(\lambda)^{r}}{r!}\\&\nonumber\quad\quad \quad\times\frac{\Gamma(\frac{1}{\rho}(\frac{\alpha}{k}+\frac{\nu r}{k}+(n-\gamma)\rho+\rho-1)}
{\Gamma(\frac{1}{\rho}(\frac{\alpha}{k}+\frac{\nu r}{k}-\gamma\rho+\rho-1)}
~s^{\frac{\alpha}{k}+\frac{\nu}{k}-\rho\gamma-1}
\end{align}
using (1.8), we obtain
\begin{align}\nonumber
P\equiv(\frac{k}{\rho})^{-\gamma}~s^{\frac{\alpha}{k}-\rho\gamma-1}
~_{n+1}\Phi_{m+1}^{k}\Bigg[\begin{matrix}\big(p_{i},\alpha_{i}\big)_{1,n},~\big(\frac{1}{\rho}(\alpha+(\rho-1)k),\frac{\nu}{\rho}\big)\quad \\
\big(q_{j},\beta_{j}\big)_{1,m},~\big(\frac{1}{\rho}(\alpha+(\rho(1-\gamma)-1)k),\frac{\nu}{\rho}\big)
\end{matrix}\Bigg|~\lambda~s^{\frac{\nu}{k}}\Bigg]\quad\quad\quad\quad\quad\Box
\end{align}
\\\textbf{\ Theorem 4.}
Let $\gamma,~ \alpha\in\mathbb{C}$ such that $\mathfrak{R}(\gamma)>0,~
\mathfrak{R}(\alpha)>1+[\mathfrak{R}(\gamma)]-\mathfrak{R}(\gamma);~ \lambda\in\mathbb{C},~ \rho>0,~\nu>0,$
then for $\Delta>-1,$ the Katugampola fractional integration $_{\rho}D_{-}^{\gamma}$ for generalized $k-$Wright
function $_{n}\Phi_{m}^{k}(z)$ is given as
\begin{align}\nonumber
\Bigg(~_{\rho}D_{-}^{\gamma}&\Bigg(\tau^{-\frac{\alpha}{k}}~_{n}\Phi_{m}^{k}
\Bigg[\begin{matrix}(p_{i},\alpha_{i})_{1,n}\\(q_{j},\beta_{j})_{1,m}
\end{matrix}\Big|~\lambda~\tau^{-\frac{\nu}{k}}\Bigg]\Bigg) \Bigg)(s)\\&
=(\frac{k}{\rho})^{-\gamma}~s^{-\rho\gamma-\frac{\alpha}{k}}
~_{n+1}\Phi_{m+1}^{k}\Bigg[\begin{matrix}\big(p_{i},\alpha_{i}\big)_{1,n},~\big(\frac{\alpha}{\rho}+k\gamma,\frac{\nu}{\rho}\big)\quad \\
\big(q_{j},\beta_{j}\big)_{1,m},~\big(\frac{\alpha}{\rho},\frac{\nu}{\rho}\big)
\end{matrix}\Bigg|~\lambda~s^{-\frac{\nu}{k}}\Bigg]\label{e4.2}
\end{align}
\\\textbf{\ Proof.} According to Lemma 1, a generalized $k-$Wright
functions in both side of the equation (4.2), exist for $s>0.$
Let $n=1+[\mathfrak{R}(\gamma)],$ then we consider that
\begin{align}\nonumber
Q\equiv\Bigg(~_{\rho}D_{-}^{\gamma}\Bigg(\tau^{-\frac{\alpha}{k}}~_{n}\Phi_{m}^{k}
\Bigg[\begin{matrix}(p_{i},\alpha_{i})_{1,n}\\(q_{j},\beta_{j})_{1,m}
\end{matrix}\Big|~\lambda~\tau^{-\frac{\nu}{k}}\Bigg]\Bigg) \Bigg)(s)
\end{align}
using (1.4), we have
\begin{align}\nonumber
Q\equiv\big(-s^{1-\rho}\frac{d}{ds}\big)^{n}\Bigg(~_{\rho}I_{-}^{n-\gamma}
\Bigg(\tau^{-\frac{\alpha}{k}}~_{n}\Phi_{m}^{k}
\Bigg[\begin{matrix}(p_{i},\alpha_{i})_{1,n}\\(q_{j},\beta_{j})_{1,m}
\end{matrix}\Big|~\lambda~\tau^{-\frac{\nu}{k}}\Bigg]\Bigg) \Bigg)(s)
\end{align}
using Theorem 2, we obtain
\begin{align}\nonumber
Q\equiv\big(-s^{1-\rho}\frac{d}{ds}\big)^{n}
(\frac{k}{\rho})^{n-\gamma}~s^{\rho(n-\gamma)-\frac{\alpha}{k}}
~_{n+1}\Phi_{m+1}^{k}\Bigg[\begin{matrix}\big(p_{i},\alpha_{i}\big)_{1,n},~\big(\frac{\alpha}{\rho}-k(n-\gamma),\frac{\nu}{\rho}\big)\quad \\
\big(q_{j},\beta_{j}\big)_{1,m},~\big(\frac{\alpha}{\rho},\frac{\nu}{\rho}\big)
\end{matrix}\Bigg|~\lambda~s^{-\frac{\nu}{k}}\Bigg]
\end{align}
using (1.12), we can write the above equation as
\begin{align}\nonumber
Q\equiv(-1)^{n}(\frac{k}{\rho})^{n-\gamma}
\sum_{r=0}^{\infty}\frac{\prod_{i=1}^{n}~\Gamma_{k}(p_{i}+\alpha_{i}r)\Gamma_{k}(\frac{\alpha}{\rho}-(n-\gamma)k+\frac{\nu}{\rho}r)}
{\prod_{j=1}^{m}~\Gamma_{k}(q_{j}+\beta_{j}r)\Gamma_{k}(\frac{\alpha}{\rho}+\frac{\nu}{\rho}r)}
\frac{(\lambda)^{r}}{r!}\big(s^{1-\rho}\frac{d}{ds}\big)^{n}\big(~s^{\rho(n-\gamma)-\frac{\alpha}{k}-\frac{\nu}{k}}\big)
\end{align}
on simplifying the above equation, we obtain
\begin{align}\nonumber
Q\equiv(-1)^{n}k^{n-\gamma}\rho^{\gamma}
\sum_{r=0}^{\infty}&\frac{\prod_{i=1}^{n}~\Gamma_{k}(p_{i}+\alpha_{i}r)\Gamma_{k}(\frac{\alpha}{\rho}-(n-\gamma)k+\frac{\nu}{\rho}r)}
{\prod_{j=1}^{m}~\Gamma_{k}(q_{j}+\beta_{j}r)\Gamma_{k}(\frac{\alpha}{\rho}+\frac{\nu}{\rho}r)}
\frac{(\lambda)^{r}}{r!}\\&\nonumber\quad\quad \quad\times
\frac{\Gamma(1+(n-\gamma)-\frac{\alpha}{\rho k}-\frac{\nu}{\rho k}r)}{\Gamma(1-\gamma-\frac{\alpha}{\rho k}-\frac{\nu}{\rho k}r)}
\big(~s^{-\rho\gamma-\frac{\alpha}{k}-\frac{\nu}{k}}\big)
\end{align}
using (1.8), we obtain
\begin{align}\nonumber
Q\equiv(-1)^{n}\rho^{\gamma}
\sum_{r=0}^{\infty}&\frac{\prod_{i=1}^{n}~\Gamma_{k}(p_{i}+\alpha_{i}r)}
{\prod_{j=1}^{m}~\Gamma_{k}(q_{j}+\beta_{j}r)\Gamma(\frac{\alpha}{\rho k}+\frac{\nu}{\rho k}r)}
\frac{(\lambda)^{r}}{r!}\\&\quad\quad \quad\times
\frac{\Gamma(\gamma-n+\frac{\alpha}{\rho k}+\frac{\nu}{\rho k}r)\Gamma(1-(\gamma-n+\frac{\alpha}{\rho k}+\frac{\nu}{\rho k}r))}
{\Gamma(1-(\gamma+\frac{\alpha}{\rho k}+\frac{\nu}{\rho k}r))}
\big(~s^{-\rho\gamma-\frac{\alpha}{k}-\frac{\nu}{k}}\big)\label{e4.3}
\end{align}
using (1.9), we have
\begin{align}\nonumber\quad\quad \quad
\Gamma(\gamma-n+\frac{\alpha}{\rho k}+\frac{\nu}{\rho k}r)&\Gamma(1-(\gamma-n+\frac{\alpha}{\rho k}+\frac{\nu}{\rho k}r))\\&\nonumber=
\frac{\pi}{\sin[(\gamma+\frac{\alpha}{\rho k}+\frac{\nu}{\rho k}r)\pi-n\pi]}\\&\nonumber=
\frac{\pi}{\sin[(\gamma+\frac{\alpha}{\rho k}+\frac{\nu}{\rho k}r)\pi]\cos(n\pi)}\\&=
\frac{(-1)^{n}\pi}{\sin[(\gamma+\frac{\alpha}{\rho k}+\frac{\nu}{\rho k}r)\pi]}\label{e4.4}
\end{align}
and
\begin{align}
\frac{1}{\Gamma(1-(\gamma+\frac{\alpha}{\rho k}+\frac{\nu}{\rho k}r))}
=\frac{\Gamma(\gamma+\frac{\alpha}{\rho k}+\frac{\nu}{\rho k}r)\sin[(\gamma+\frac{\alpha}{\rho k}+\frac{\nu}{\rho k}r)\pi]}{\pi}\label{e4.5}
\end{align}
Substituting (4.4) and (4.5) in (4.3), and finally using (1.8), we obtain
\begin{align}\nonumber
Q\equiv(\frac{k}{\rho})^{-\gamma}~s^{-\rho\gamma-\frac{\alpha}{k}}
~_{n+1}\Phi_{m+1}^{k}\Bigg[\begin{matrix}\big(p_{i},\alpha_{i}\big)_{1,n},~\big(\frac{\alpha}{\rho}+k\gamma,\frac{\nu}{\rho}\big)\quad \\
\big(q_{j},\beta_{j}\big)_{1,m},~\big(\frac{\alpha}{\rho},\frac{\nu}{\rho}\big)
\end{matrix}\Bigg|~\lambda~s^{-\frac{\nu}{k}}\Bigg]
\end{align}
\\\textbf{\ Remark 2.} If $\rho=1,$ then
\par Theorems 1, 2, 3 and 4, are reduced to Theorems 2, 3, 4 and 5, respectively (see[10]).

\end{document}